\documentclass[12pt]{article}

\usepackage{amssymb}
\usepackage{amsmath}
\usepackage{graphicx}
\usepackage{hyperref}

\numberwithin{equation}{section}

\begin{document}

\newcommand{\bd}{\begin{displaymath}}
\newcommand{\ed}{\end{displaymath}}
\newcommand{\ds}{\displaystyle}
\newcommand{\bp}{\underline{\bf Proof}:\ }
\newcommand{\ep}{{\hfill $\Box$}\\ }
\newcommand{\be}{\begin{equation}}
\newcommand{\ee}{\end{equation}}
\newcommand{\ba}{\begin{array}}
\newcommand{\ea}{\end{array}}
\newcommand{\bea}{\begin{eqnarray}}
\newcommand{\eea}{\end{eqnarray}}
\newcommand{\nt}{\noindent}

\newtheorem{0}{DEFINITION}[section]
\newtheorem{1}{LEMMA}[section]
\newtheorem{2}{THEOREM}[section]
\newtheorem{3}{COROLLARY}[section]
\newtheorem{4}{PROPOSITION}[section]
\newtheorem{5}{REMARK}[section]
\newtheorem{6}{EXAMPLE}[section]
\newtheorem{7}{ALGORITHM}[section]
\newtheorem{8}{CONJECTURE}[section]
\newtheorem{9}{QUESTION}[section]

\title{Can a Higher Order Markov Chain Be Treated as a First Order Markov Chain?}
\author{
Jianhong Xu\thanks{School of Mathematical and Statistical Sciences, Southern Illinois University Carbondale, Carbondale, IL 62901, USA. Email: \texttt{jhxu@siu.edu}} 
% \and  \thanks{}
% \and  \thanks{}
}

\maketitle

\begin{abstract}  
It is well known that any higher order Markov chain can be associated with a first order Markov chain. In this primarily expository article, we present the first fairly comprehensive analysis of the relationship between higher order and first order Markov chains, together with illustrative examples. Our main objective is to address the central question as posed in the title.
\end{abstract}

\nt {\bf Keywords}: 
higher order Markov chain, transition tensor, $k$-step transition probability, ever-reaching probability, ergodicity, regularity, mean first passage time, limiting probability distribution, stationary distribution, recurrence, transience, classification of states

\nt {\bf AMS Subject Classification}: 60A05, 60J10, 60J99, 15A69, 15B51

\section{Introduction}
\label{intro}
\setcounter{equation}{0}

A first order or simple Markov chain, or just Markov chain as is often called in literature, is a widely used stochastic model in which the future state of a system is determined by its current state alone via some probabilistic rules, not by any of its past states. In many real-world problems, however, it often arises that the future state of a system depends on not only its current state
but also a number of its past states following certain probabilistic rules. A higher order or multiple Markov chain is exactly a stochastic model for such systems and has been adopted in many applied areas \cite{Bae, Bur, Fle, Ho, Isl, Kwa, Lan, LLZ, M, San, Xio, YJK}. It is a generalization of the notion of first order Markov chains.

Formally, first order and higher order Markov chains can be defined in a unified way as follows. Let $m \ge 2$. Then, an $(m-1)$th order Markov chain is a stochastic process $\{X_t : t=1, 2, \ldots\}$, where $t$ can be regarded as time and $X_t$ is a random variable representing the state of a system at $t$ while taking values in the state space $S=\{1, 2, \ldots, n\}$ with $n \ge 2$, such that
\be
\label{markov}
\begin{array}{l}
\Pr(X_{t+1}=i_1 | X_t=i_2,\ldots,X_{t-m+2}=i_m,\ldots,X_1=i_{t+1})\\
=\Pr(X_{t+1}=i_1 | X_t=i_2,\ldots,X_{t-m+2}=i_m)
\end{array}
\ee
for any $t \ge m-1$ and $i_1,i_2,\ldots,i_m,\ldots,i_{t+1} \in S$. In particular, this chain is called higher order when $m \ge 3$. Incidentally, in this definition, $X_{t+1}$ can be interpreted as the future state of the system, $X_t$ the current state, and $X_{t-1}$ through $X_{t-m+2}$ the $m-2$ past states. By reason of the latter, an $(m-1)$th order Markov chain is also said to have memory $m-2$ and, especially, to be memoryless if it is of first order with $m=2$.

In addition, throughout this article, a Markov chain is always assumed to be homogeneous, i.e., the probability in (\ref{markov}) depends only on $i_1, i_2, \ldots, i_m$ but is independent of $t$. Accordingly, we can write there 
$$\Pr(X_{t+1}=i_1 | X_t=i_2,\ldots,X_{t-m+2}=i_m)=p_{i_1i_2\ldots i_m},$$
which is called the transition probability from states $(i_2,\ldots,i_m)$ to $i_1$. These probabilities form an $m$th order, $n$ dimensional tensor, which is denoted by ${\cal P}=[p_{i_1i_2\ldots i_m}]$ and called the transition tensor. Such a $\cal P$ is stochastic in the sense that $0 \le p_{i_1i_2\ldots i_m} \le 1$ for all $i_1, i_2,\ldots,i_m \in S$ and $\displaystyle \sum_{i_1 \in S} p_{i_1i_2\ldots i_m}=1$ for all $i_2, \ldots, i_m \in S$. For a first order Markov chain, $\cal P$ is reduced to an $n \times n$ (column stochastic) transition matrix, denoted as $P$ or $Q$ in this article.

For brevity, in what follows, we shall refer to a Markov chain of any order simply as a chain.

From the foregoing unified definition, it is natural to anticipate any result for higher order chains to extend to the first order case as well. The converse, nevertheless, is not so clear. Specifically, a key question is whether a problem regarding higher order chains can be solved by way of first order chains. As far as we can see, it has been a prevailing perception that the answer to this question is to a large extent affirmative. For one thing, it is well known that a higher order chain can be associated with a first order chain \cite{Doo, Hun}. Detail in the respect will be given in the next section.

The main thrust of this article is to shed light on the above key question. The importance of this question lies in the fact that it significantly impacts the methodology of examining higher order chains. Even though such chains have lately been the focus of many studies, see \cite{CZ, CPZ, G, GLY, HQ, LZ, LN, WC} and the references therein, this question remains open to a great degree. In light of this situation, we shall provide, for the first time, a rather comprehensive answer by synthesizing relevant recent results from \cite{HWX, HX24a, HX24b, HX25}. We shall also give necessary examples to illustrate the results. In particular, we shall highlight the primary differences between higher order and first order chains. These differences, to the best of our knowledge, have yet to be presented in such a collective and systematic way.

To help verify the examples in this article or related works in \cite{HWX, HX24a, HX24b, HX25}, the reader may use the HOMC MATLAB package \cite{X25}, which is publicly available at 

\begin{center}
\url{https://neumann.math.siu.edu/homc}.
\end{center}

\section{Reduced First Order Chains}
\label{reduc}
\setcounter{equation}{0}

Let $m \ge 3$. Suppose that $X=\{X_t : t=1, 2, \ldots\}$ is an $(m-1)$th order chain on state space $S=\{1, 2, \ldots, n\}$ whose transition tensor is ${\cal P}=[p_{i_1i_2\ldots i_m}]$. It is well known \cite{Doo, Hun} that this higher order chain $X$ can be associated with a first order chain as follows.

First, define the set $T$ of multi-indices of length $m-1$, ordered by linear indexing \cite{MSL}, as $$T=\{i_1i_2\ldots i_{m-1} : i_1, i_2, \ldots, i_{m-1} \in S\}.$$
Note that the cardinality of $T$ is $N=n^{m-1}$. Next, for $t \ge m-1$, introduce vector-valued random variables 
$$Y_t=[X_t \ \ X_{t-1} \ \ \ldots \ \ X_{t-m+2}]^T.$$ 
Let us denote $Y_t=i_1i_2\ldots i_{m-1}$ when $X_t=i_1, X_{t-1}=i_2, \ldots, X_{t-m+2}=i_{m-1}$. Then, the vector-valued stochastic process $Y=\{Y_t : t=m-1, m, \ldots\}$ is a first order chain on state space $T$, which we shall refer to as the reduced first order chain obtained from $X$.

Denote the $N \times N$ transition matrix of $Y$ by $Q$. It is more convenient to write 
$$Q=[q_{i_1i_2\ldots i_{m-1},j_2j_3\ldots j_m}]$$ 
in multi-index form with $i_1i_2\ldots i_{m-1}, j_2j_3\ldots j_m \in T$. Then, the entries of $Q$ can be expressed as, see \cite{HX25},
$$q_{i_1i_2\ldots i_{m-1}, j_2j_3\ldots j_m}=\left\{\begin{array}{cl}
p_{i_1i_2\ldots i_{m-1}j_m}, & ~i_\ell=j_\ell, ~\ell=2, 3, \ldots, m-1;\\
0, & ~{\rm otherwise.}
\end{array}\right.$$
In addition, $p_{i_1i_2\ldots i_m}$ is located in $Q$ on row 
$$i_1+n(i_2-1)+\ldots +n^{m-2}(i_{m-1}-1)$$
and column 
$$i_2+n(i_3-1)+\ldots +n^{m-2}(i_m-1).$$
For a third order chain with two states, i.e., $m=4$ and $n=2$, for example, we have 
$$T=\{111, 211, 121, 221, 112, 212, 122, 222\}$$ and 
$$Q=\left[\ba{cccccccc}
p_{1111} & 0 & 0 & 0 & p_{1112} & 0 & 0 & 0\\
p_{2111} & 0 & 0 & 0 & p_{2112} & 0 & 0 & 0\\
0 & p_{1211} & 0 & 0 & 0 & p_{1212} & 0 & 0\\
0 & p_{2211} & 0 & 0 & 0 & p_{2212} & 0 & 0\\
0 & 0 & p_{1121} & 0 & 0 & 0 & p_{1122} & 0\\
0 & 0 & p_{2121} & 0 & 0 & 0 & p_{2122} & 0\\
0 & 0 & 0 & p_{1221} & 0 & 0 & 0 & p_{1222}\\
0 & 0 & 0 & p_{2221} & 0 & 0 & 0 & p_{2222}
\ea\right].$$
In particular, $Q$ coincides with the transition matrix of $X$ when this chain is first order, i.e., in the special case $m=2$.

Speaking of the connection between a higher order chain and its reduced first order chain, we mention here the opinion of Joseph Doob, who has left a profound, lasting impact on the field of stochastic processes. In the seminal monograph \cite{Doo}, he stated that ``the generalization (from simple Markov processes to multiple ones) is not very significant, because the (vector) process with random variables $\{\hat{x}_n\}$, ${\hat x}_n=(x_n, \ldots, x_{n+v-1})$ (i.e., $Y_t$ in our notation) has the Markov property ... Thus multiple Markov processes can be reduced to simple ones at the small expense of going to vector-valued random variables.'' This appears to suggest that problems regarding a higher chain can be solved through its reduced first order chain.

From what we can gather in literature, different voices have been few and far between. The most notable argument can be found in the classic treatise \cite{Ios} by Marius Iosifescu, in which he pointed out first that ``In this way with any multiple Markov chain one can associate a simple Markov chain. This is the origin of the rooted prejudice that the study of multiple Markov chains reduces to that of simple ones.'' and then cited the following example: State $i$ can be recurrent in a second order chain but $ij \in T$ is not recurrent in the reduced first order chain for any $j$. This example, albeit useful, falls short to adequately address why it is a ``prejudice'' to investigate higher order chains by reducing them to first order chains, i.e., to fully address the key question as raised in the introductory section. Having laid out the notion of reduced first order chains, we can now rephrase this key questions as:
\begin{9}
\label{ques}
Can a higher order chain be studied via its reduced first order chain?
\end{9}
This formal statement clarifies what the first order chain is in the previously stated approach ``by way of first order chains''. It is also easier for us to refer to this question in the rest of this article.

\section{$k$-Step Transition and Ever-Reaching Probabilities}
\label{prob}
\setcounter{equation}{0}

Let us begin here with some background material.

For $m, n \ge 2$, let ${\cal A}=[a_{i_1i_2\ldots i_m}]$ and ${\cal B}=[b_{i_1i_2\ldots i_m}]$ be both $m$th order, $n$ dimensional tensors. The tensor product\footnote{This product was initially denoted by ${\cal B} \boxtimes {\cal A}$ in \cite{HX24a}. The current notation is consistent with the special case when $\cal A$ and $\cal B$ are matrices.} ${\cal A} \boxtimes {\cal B}$ is introduced in \cite{HX24a, HX25} as an $m$th order, $n$ dimensional tensor ${\cal C}=[c_{i_1i_2\ldots i_m}]$ such that
$$c_{i_1i_2\ldots i_m}=\sum_{j=1}^n a_{i_1ji_2\ldots i_{m-1}}b_{ji_2\ldots i_m}$$
for any $1 \le i_1, i_2, \ldots, i_m \le n$. When $m=2$, in particular, ${\cal A}\boxtimes {\cal B}$ reduces to a matrix multiplication. Note, however, that the $\boxtimes$ product is in general not associative whenever $m \ge 3$.

For $k=2, 3, \ldots$, the $k$th power of the $m$th order, $n$ dimensional tensor $\cal A$ is defined recursively by $${\cal A}^{k+1}={\cal A}^k \boxtimes {\cal A}, ~k=1, 2, \ldots,$$ with ${\cal A}^1=\cal A$. By convention, ${\cal A}^0={\cal I}=[\delta_{i_1i_2\ldots i_m}]$, the $m$th order, $n$ dimensional identity tensor whose entries satisfy 
$$\delta_{i_1i_2\ldots i_m}=\left\{\begin{array}{cl}
1, & i_1=i_2;\\
0, & {\rm otherwise.}
\end{array}\right.$$
It is easy to see that ${\cal I}\boxtimes {\cal A}=\cal A$, but in general ${\cal A}\boxtimes {\cal I} \ne \cal A$.

Besides the above, the mode-$k$ matricization of the $m$th order, $n$ dimensional tensor $\cal A$ is defined to be an $n$ by $N$, with $N=n^{m-1}$, matrix that has the mode-$k$ fibers ${\cal A}(i_1, \ldots, i_{k-1}, :, i_{k+1}, \ldots, i_m)$ as columns. These fibers are arranged in the linear indexing order of $i_1\ldots i_{k-1}i_{k+1}\ldots i_m$ \cite{KB}. Especially, the mode-$1$ matricization of $\cal A$ is formed by arranging ${\cal A}(:,:,i_3, \ldots, i_m)$, i.e., the frontal slices of $\cal A$, side by side via the linear indexing order of $i_3\ldots i_m$.

As shown in \cite{HWX, HX24a, HX25}, the $\boxtimes$ product plays a crucial role in the study of higher order chains. It can be used to formulate various important quantities such as $k$-step transition probabilities and ever-reaching probabilities, which we shall describe next. The mode-$1$ matricization, meanwhile, is useful when dealing with the limiting probability distribution problem.

Back to the $(m-1)$th order chain $X$ with state space $S=\{1, 2, \ldots, n\}$ and transition tensor ${\cal P}=[p_{i_1i_2\ldots i_m}]$. Denote ${\cal P}^k=[p^{(k)}_{i_1i_2\ldots i_m}]$. Then, $p^{(k)}_{i_1i_2\ldots i_m}$ is the $k$-step transition probability from states $(i_2, \ldots, i_m)$ to $i_1$ \cite{HWX}, i.e.,
$$p^{(k)}_{i_1i_2\ldots i_m}=\Pr(X_{t+k}=i_1 | X_t=i_2, \ldots, X_{t-m+2}=i_m).$$
Accordingly, ${\cal P}^k$ is called the $k$-step transition tensor.

Continuing, let 
\be
\label{eta}
\eta_{i_1i_2\ldots i_m}=\min\{j \ge 1 : X_{m+j-1}=i_1 | X_{m-1}=i_2, \ldots, X_1=i_m\}
\ee
be the random variable of the first passage time from states $(i_2, \ldots, i_m)$ to $i_1$. Denote the probability for this passage to happen at $j=k$ by 
$$f^{[k]}_{i_1i_2\ldots i_m}=\Pr(\eta_{i_1i_2\ldots i_m}=k).$$
In the special case when $k=1$, $f^{[1]}_{i_1i_2\ldots i_m}=p_{i_1i_2\ldots i_m}$. Denote ${\cal F}^{[k]}=[f^{[k]}_{i_1i_2\ldots i_m}]$, which is called the $k$-step first passage time probability tensor. It is known \cite{HWX} that  
\be
\label{f}
f^{[k+1]}_{i_1i_2\ldots i_m}=\sum_{j \in S, j \ne i_1}f^{[k]}_{i_1ji_2\ldots i_{m-1}}p_{ji_2\ldots i_m}, ~k=1, 2, \ldots,
\ee
To express (\ref{f}) via the $\boxtimes$ product, for any $m$th order, $n$ dimensional tensor ${\cal A}=[a_{i_1i_2\ldots i_m}]$, let us define the diagonal tensor ${\cal A}_d=[a^{(d)}_{i_1i_2\ldots i_m}]$ satisfying $$a^{(d)}_{i_1i_2\ldots i_m}=\left\{\begin{array}{cl}
a_{i_1i_1i_3\ldots i_m}, & i_1=i_2;\\
0, & {\rm otherwise}.
\end{array}\right.$$
Clearly, (\ref{f}) can now be written simply as 
$${\cal F}^{[k+1]}=({\cal F}^{[k]}-{\cal F}^{[k]}_d)\boxtimes \cal P.$$
The probability of ever reaching $i_1$ from $(i_2, \ldots, i_m)$ is given by  
$$f_{i_1i_2\ldots i_m}=\Pr(\eta_{i_1i_2\ldots i_m}<\infty)=\sum_{k=1}^\infty f^{[k]}_{i_1i_2\ldots i_m},$$
or in tensor form, the ever-reaching probability tensor ${\cal F}=[f_{i_1i_2\ldots i_m}]$ can be expressed as  
$${\cal F}=\sum_{k=1}^\infty {\cal F}^{[k]}.$$

These results regarding $k$-step transition probabilities and ever-reaching probabilities remain valid for first order chains. In the first order case, they can also be written as $P^k=[p^{(k)}_{ij}]$, the $k$-step transition matrix, and $F=[f_{ij}]$, the ever-reaching probability matrix, respectively \cite{Hun, KS}.

In view of Question \ref{ques}, we now raise the question of whether the $k$-step transition tensor or ever-reaching probability tensor of a higher order chain can be obtained from their respective counterpart of the reduced first order chain. We discuss the $k$-step transition probabilities here. The ever-reaching probabilities will be considered in the next section after the introduction of the notion of ergodicity. 

Let us take, for example, a third order chain with transition tensor ${\cal P}=[p_{i_1i_2i_3i_4}]$ and state space $S=\{1, 2, \ldots, n\}$. Suppose that $Q=[q_{i_1i_2i_3, j_1j_2j_3}]$ is the transition matrix of the reduced first order chain, whose state space is $T=\{i_1i_2i_3 : i_1, i_2, i_3 \in S\}$. 
\begin{2}
\label{pkthm}
 Under the above settings, 
 $$p^{(2)}_{i_1i_2i_3i_4}=\sum_{j_1 \in S}q^{(2)}_{i_1j_1i_2, i_2i_3i_4}$$
 and, for $k \ge 3$, 
 \be
 \label{pk}
 p^{(k)}_{i_1i_2i_3i_4}=\sum_{j_2 \in S}\sum_{j_1 \in S}q^{(k)}_{i_1j_1j_2, i_2i_3i_4}.
 \ee
\end{2}
\bp
From $$q^{(2)}_{i_1i_2i_3, j_1j_2j_3}=\sum_{k_1k_2k_3 \in T}q_{i_1i_2i_3, k_1k_2k_3}q_{k_1k_2k_3, j_1j_2j_3},$$
we notice that the terms in this summation are all zero except when $k_1=i_2$, $k_2=i_3=j_1$, and $k_3=j_2$, i.e.,
\be
\label{q2}
q^{(2)}_{i_1i_2i_3, i_3j_2j_3}=p_{i_1i_2i_3j_2}p_{i_2i_3j_2j_3}.
\ee
It follows that 
$$p^{(2)}_{i_1i_2i_3i_4}=\sum_{j_1 \in S}p_{i_1j_1i_2i_3}p_{j_1i_2i_3i_4}=\sum_{j_1 \in S}q^{(2)}_{i_1j_1i_2, i_2i_3i_4}.$$

Similarly, we observe that in  
$$q^{(3)}_{i_1i_2i_3, j_1j_2j_3}=\sum_{k_1k_2k_3 \in T}q^{(2)}_{i_1i_2i_3, k_1k_2k_3}q_{k_1k_2k_3, j_1j_2j_3},$$
all the terms are zero except when $k_1=i_3$, $k_2=j_1$, and $k_3=j_2$, i.e., 
$$q^{(3)}_{i_1i_2i_3, j_1j_2j_3}=p_{i_1i_2i_3j_1}p_{i_2i_3j_1j_2}p_{i_3j_1j_2j_3}.$$
Accordingly, 
$$p^{(3)}_{i_1i_2i_3i_4}=\sum_{j_2 \in S}p^{(2)}_{i_1j_2i_2i_3}p_{j_2i_2i_3i_4}=\sum_{j_2 \in S}\sum_{j_1 \in S}p_{i_1j_1j_2i_2}p_{j_1j_2i_2i_3}p_{j_2i_2i_3i_4}=\sum_{j_2 \in S}\sum_{j_1 \in S}q^{(3)}_{i_1j_1j_2, i_2i_3i_4}.$$

Next, suppose that (\ref{pk}) holds for $k=3, 4, \ldots, \ell$. Then,
$$\begin{array}{rcl}
p^{(\ell +1)}_{i_1i_2i_3i_4} & = & \ds \sum_{k_1 \in S}p^{(\ell)}_{i_1k_1i_2i_3}p_{k_1i_2i_3i_4}\\
 & = & \ds \sum_{j_2 \in S}\sum_{j_1 \in S}\sum_{k_1 \in S}q^{(\ell)}_{i_1j_1j_2, k_1i_2i_3}q_{k_1i_2i_3, i_2i_3i_4}\\
 & = & \ds \sum_{j_2 \in S}\sum_{j_1 \in S}\sum_{k_1k_2k_3 \in T}q^{(\ell)}_{i_1j_1j_2,k_1k_2k_3}q_{k_1k_2k_3, i_2i_3i_4},
\end{array}$$
where the last equality is due to $q_{k_1k_2k_3,i_2i_3i_4}=0$ if $k_2 \ne i_2$ or $k_3 \ne i_3$. Hence, 
$$p^{(\ell +1)}_{i_1i_2i_3i_4}=\sum_{j_2 \in S}\sum_{j_1 \in S}q^{(\ell +1)}_{i_1j_1j_2, i_2i_3i_4}.$$
The conclusion in (\ref{pkthm}) follows now by induction.
\ep

Theorem \ref{pkthm} provides a connection between the $k$-step transition probabilities of a third order chain and those of its reduced first order chain. This result can be extended to a general higher order chain as well. We, however, shall skip this since the $k$-step transition probabilities of such a chain can be obtained directly from ${\cal P}^k$, whereas the expression for $q^{(k)}_{i_1i_2\ldots i_{m-1}, j_1j_2\ldots j_{m-1}}$ in terms of the entries of $\cal P$ becomes rather complicated and multiple summations are needed to recover $p^{(k)}_{i_1i_2\ldots i_m}$ as $m$ and $k$ increase.

To end this section, let us give one example of the mode-$1$ matricization of the transition tensor $\cal P$, denoted by $P$. Take a third order chain with two states just as the one in the previous section, i.e., $m=4$ and $n=2$. Then,  
$$P=\left[\ba{cc|cc|cc|cc}
p_{1111} & p_{1211} & p_{1121} & p_{1221} & p_{1112} & p_{1212} & p_{1122} & p_{1222}\\
p_{2111} & p_{2211} & p_{2121} & p_{2221} & p_{2112} & p_{2212} & p_{2122} & p_{2222}
\ea\right].$$
Clearly, $P$ coincides with the transition matrix when the chain is first order.

\section{Irreducibility, Ergodicity, and Regularity}
\label{irre}
\setcounter{equation}{0}

Let us begin with several well-known concepts for first order chains \cite{Hun, KS}. Suppose that $P=[p_{ij}]$ is the transition matrix of a first order chain on state space $S$. Then:
\begin{itemize}
\item{The chain is irreducible if for any $\emptyset \ne K \subsetneq S$, there exists $p_{ij}>0$ for some $i \in K$ and some $j \in K^c$. In matrix theory, its transition matrix $P$ is also said to be irreducible. The chain is said to be reducible if it is not irreducible.}
\item{The chain is ergodic if for any $i, j \in S$, there exists $k \ge 1$, which may depend on $i$ and $j$, such that $p^{(k)}_{ij}>0$.}
\item{The chain is regular if there exists some $k \ge 1$ such that $p^{(k)}_{ij}>0$ for all $i, j \in S$, i.e., $P^k > 0$. In matrix theory, such $P$ is said to be primitive.}
\end{itemize}
For a first order chain, it is well known \cite{HJ} that irreducibility is equivalent to ergodicity.

The above definitions can be extended as follows to a higher order chain on state space $S$, whose transition tensor is ${\cal P}=[p_{i_1i_2\ldots i_m}]$.
\begin{itemize}
\item{The chain is said to be irreducible if for each $\emptyset \ne K \subsetneq S$, there exists $p_{i_1i_2\ldots i_m}>0$ for some $i_1 \in K$ and some $i_2,\ldots,i_m \in K^c$ \cite{LN}.}
\item{The chain is said to be ergodic if for any $i_1, i_2, \ldots, i_m \in S$, there exists $k \ge 1$, which may depend on $i_1, i_2, \ldots, i_m$, such that $p^{(k)}_{i_1i_2\ldots i_m}>0$ \cite{HWX}.}
\item{The chain is said to be regular if there exists $k \ge 1$ such that for any $i_1, i_2, \ldots, i_m \in S$, $p^{(k)}_{i_1i_2\ldots i_m}>0$, i.e., ${\cal P}^k > 0$ \cite{HX25}.}
\end{itemize}
Unlike the first order case, for a higher order chain, ergodicity is a stronger condition than irreducibility \cite{HWX}, i.e., an ergodic higher order chain must be irreducible, but the converse is not true in general. Consider, for example, a second order chain with three states and transition tensor 
$${\cal P}(:,:,1)=\left[\begin{array}{ccc}
0 & 0 & 0\\
1 & 0 & 0\\
0 & 1 & 1
\end{array}\right], ~{\cal P}(:,:,2)=\left[\begin{array}{ccc}
0 & 0 & 0\\
0 & 0 & 0\\
1 & 1 & 1
\end{array}\right], ~{\cal P}(:,:,3)=\left[\begin{array}{ccc}
0 & 0 & 1\\
0 & 0 & 0\\
1 & 1 & 0
\end{array}\right].$$
This chain is irreducible, but not ergodic since $p^{(k)}_{2i_2i_3}=0$ for all $1 \le i_2, i_3 \le 3$ and $k \ge 2$.

In order to guarantee the mean first passage times to be well defined and finite, for example, ergodicity is an important condition no matter the order of a chain \cite{HX24a, KS}. Following Question \ref{ques}, it is natural to ask whether the mean first passage times of a higher order chain can be obtained from those of its reduced first order chain. First and foremost, however, it makes sense to ask whether ergodicity of a higher order chain carries over to its reduced first order chain. The answer to the latter question, unfortunately, turns out to be negative as shown by the following example.

Let ${\cal P}=[p_{i_1i_2i_3}]$ be the transition tensor of a second order chain of three states such that all $p_{i_1i_2i_3}$ are positive except $p_{311}=p_{312}=p_{313}=0$. Let us simply take  
$${\cal P}(:,:,1)={\cal P}(:,:,2)={\cal P}(:,:,3)=\left[\begin{array}{ccc}
1/2 & 1/3 & 1/3\\
1/2 & 1/3 & 1/3\\
0 & 1/3 & 1/3
\end{array}\right].$$
Since ${\cal P}^2>0$, this chain is regular and thus is also ergodic. Its reduced first order chain has transition matrix 
$$Q=\left[\begin{array}{ccccccccc} 
1/2 & 0 & 0 & 1/2 & 0 & 0 & 1/2 & 0 & 0 \\
1/2 & 0 & 0 & 1/2 & 0 & 0 & 1/2 & 0 & 0 \\
0 & 0 & 0 & 0 & 0 & 0 & 0 & 0 & 0 \\
0 & 1/3 & 0 & 0 & 1/3 & 0 & 0 & 1/3 & 0 \\
0 & 1/3 & 0 & 0 & 1/3 & 0 & 0 & 1/3 & 0 \\
0 & 1/3 & 0 & 0 & 1/3 & 0 & 0 & 1/3 & 0 \\
0 & 0 & 1/3 & 0 & 0 & 1/3 & 0 & 0 & 1/3 \\
0 & 0 & 1/3 & 0 & 0 & 1/3 & 0 & 0 & 1/3 \\
0 & 0 & 1/3 & 0 & 0 & 1/3 & 0 & 0 & 1/3
\end{array}\right],$$
which is clearly not even irreducible due to its third row of all zeros. Hence, the reduced first order chain is not ergodic.

The above example also justifies that, like ergodicity, regularity does not carry over either from a higher order chain to its reduced first order chain.

What this example demonstrates is in fact the first obstacle to studying a higher order chain through its reduced first order chain. We shall demonstrate later in Section \ref{mfpt} that other difficulties may arise even if ergodicity or regularity happens to carry over to the reduced first order chain. The issue of mean first passage times will be addressed there too.

\begin{figure}[h!]
\begin{center}
\includegraphics[trim= 0cm 3.5cm 0cm 0cm,clip,scale=.5]{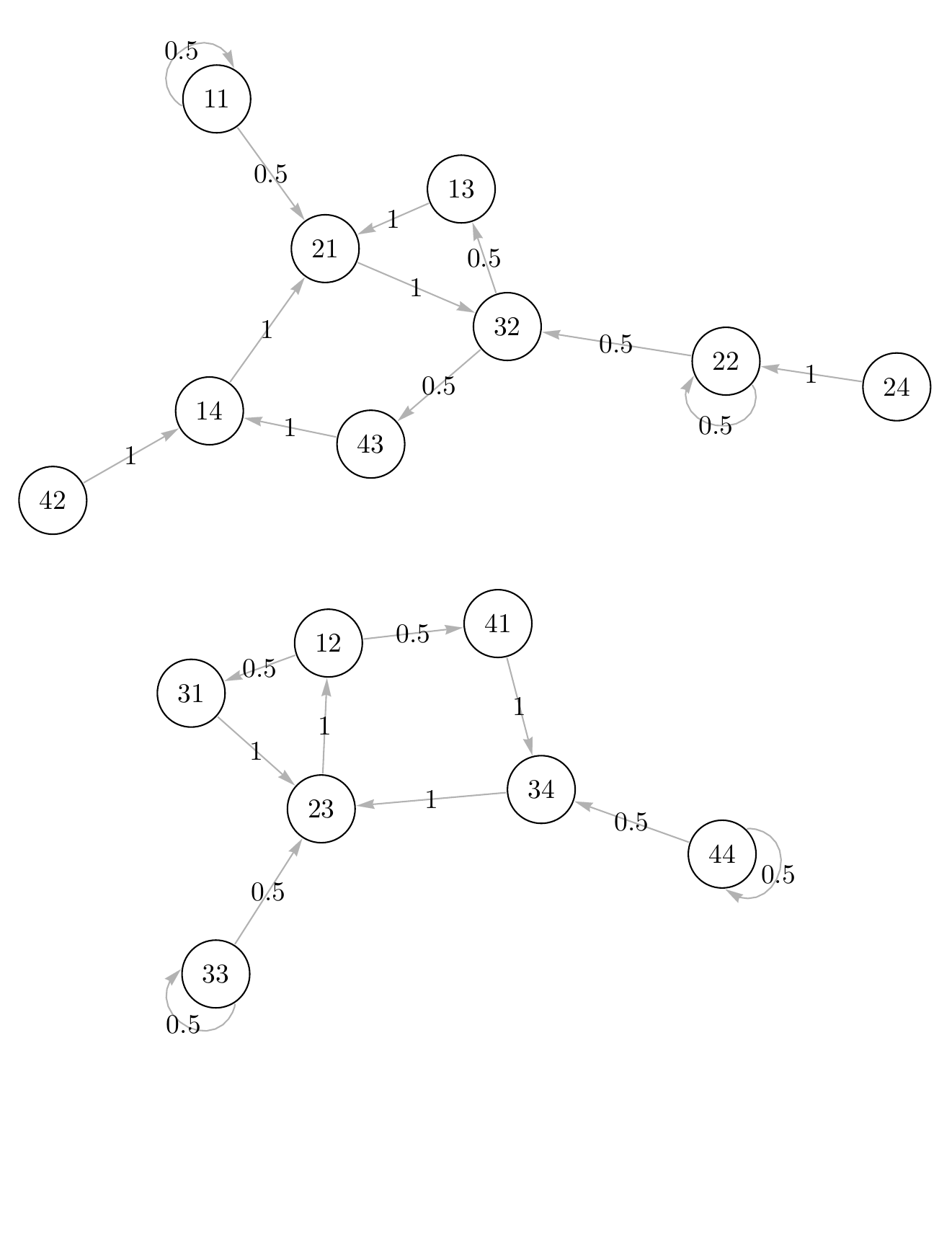}
\caption{Digraph of Reduced First Order Chain}
\label{di}
\end{center}
\end{figure}

Returning to the question raised in the preceding section concerning the relationship between the ever-reaching probabilities of a higher order chain and those of its reduced order chain, we consider below a second order chain with four states and transition tensor  
$${\cal P}(:,:,1)=\left[\ba{cccc}
1/2 & 0 & 0 & 0\\
1/2 & 0 & 1 & 0\\
0 & 1 & 0 & 1\\
0 & 0 & 0 & 0
\ea\right], ~{\cal P}(:,:,2)=\left[\ba{cccc}
0 & 0 & 1/2 & 1\\
0 & 1/2 & 0 & 0\\
1/2 & 1/2 & 0 & 0\\
1/2 & 0 & 1/2 & 0
\ea\right],$$
$${\cal P}(:,:,3)=\left[\ba{cccc}
0 & 1 & 0 & 1\\
1 & 0 & 1/2 & 0\\
0 & 0 & 1/2 & 0\\
0 & 0 & 0 & 0
\ea\right], ~{\cal P}(:,:,4)=\left[\ba{cccc}
0 & 0 & 0 & 0\\
1 & 1 & 1 & 0\\
0 & 0 & 0 & 1/2\\
0 & 0 & 0 & 1/2
\ea\right].$$
It is easy to verify ${\cal P}^{10} > 0$. Clearly, this chain is regular and, consequently, ergodic. It is shown in \cite{HX24a} that for such an ergodic chain, its ever-reaching probabilities $f_{i_1i_2i_3}=1$ for all $1 \le i_1, i_2, i_3 \le 4$.  The associated reduced first order chain, however, exhibits a noticeably different dynamics as illustrated by the digraph in Figure \ref{di} that represents its transition matrix $Q$. Starting, for example, from state $12$, the reduced first order chain has zero probability of reaching state $11$. This, nevertheless, does not lead to the conclusion that state $1$ cannot be reached from states $(1, 2)$ in the second order chain since $f_{112}=1$. This illustrates that the answer to Question \ref{ques} is negative  when dealing with the problem of ever-reaching probabilities.

It is straightforward, though a bit tedious, to show that the ever-reaching probabilities of a higher order chain cannot be determined from those of its reduced first order chain in a manner as neat as Theorem \ref{pkthm}. Consider the same settings as Theorem \ref{pkthm}, i.e., $m=4$ and $n=2$. Let $G^{[2]}=[g^{[2]}_{i_1i_2i_3, j_1j_2j_3}]$ be the $2$-step first passage time probability matrix, in multi-index notation, of the reduced first order chain. Then, 
$$g^{[2]}_{i_1i_2i_3, j_1j_2j_3}=p_{i_1i_2i_3j_2}p_{i_2i_3j_2j_3}, ~i_1i_2i_3, j_1j_2j_3 \in T,$$
provided $j_1=i_3$ and either $i_1 \ne i_2$ or $i_2 \ne i_3$ or $i_3 \ne j_2$, which is much more complicated than (\ref{q2}). This also indicates that the answer to Question \ref{ques} is negative as far as ever-reaching probabilities are concerned.

\section{Classification of States}
\label{clas}
\setcounter{equation}{0}

Given a first order chain with transition matrix $P=[p_{ij}]$ and ever-reaching probability matrix $F=[f_{ij}]$, it is well known \cite{Hun} that:
\begin{itemize}
\item{A state $i$ is recurrent if $f_{ii}=1$. A state $i$ is current iff $\displaystyle \sum_{k=1}^\infty p^{(k)}_{ii}=\infty.$}
\item{A state $i$ is transient if $f_{ii}<1$. A state $i$ is transient iff $\displaystyle \sum_{k=1}^\infty p^{(k)}_{ii}<\infty.$}
\item{A state $i$ is absorbing if $p_{ii}=1$.}
\end{itemize}

The example that is cited by \cite{Ios}, see the introductory section, indicates that the states of a higher order chain cannot be classified according to the states of its reduced first order chain. A more complete classification scheme for states of higher order chains is recently proposed in \cite{HX24b}. Here is a brief description of this scheme. Consider a higher order chain whose state space is $S$, transition tensor is ${\cal P}=[p_{i_1i_2\ldots i_m}]$, and ever-reaching probability tensor is ${\cal F}=[f_{i_1i_2\ldots i_m}]$. Then,
\begin{itemize}
\item{A state $i$ is said to be recurrent if $f_{iii_3\ldots i_m}=1$ for all $i_3, \ldots, i_m \in S$. If state $i$ is recurrent, then $\displaystyle \sum_{k=1}^\infty p^{(k)}_{iii_3\ldots i_m}=\infty$ for all $i_3, \ldots, i_m \in S$.}
\item{A state $i$ is said to be transient if $f_{iii_3\ldots i_m}<1$ for some $i_3, \ldots, i_m \in S$. Especially, a state $i$ is said to be fully transient if $f_{iii_3\ldots i_m}<1$ for any $i_3, \ldots, i_m \in S$. When state $i$ is fully transient, $\displaystyle \sum_{k=1}^\infty p^{(k)}_{iii_3\ldots i_m}<\infty$ for all $i_3, \ldots, i_m \in S$.}
\item{A state $i$ is said to be absorbing if $p_{iii_3\ldots i_m}=1$ for all $i_3, \ldots, i_m \in S$.}
\end{itemize}
These results from \cite{HX24b} seem to be parallel to their first order counterparts. We, however, point out some differences between the higher order and first order cases. Specifically, for a higher order chain: 
\begin{itemize}
\item{There may exist the scenario where a state $i$ is transient but not fully transient, i.e., $f_{iii_3\ldots i_m}$ is smaller than $1$ for some $i_3, \ldots, i_m \in S$ but is equal to $1$ for other $i_3, \ldots, i_m \in S$. This is indeed at the center of the key differences between the higher order and first order chains when it comes to classification of states.}
\item{The condition of whether $\displaystyle \sum_{k=1}^\infty p^{(k)}_{iii_3\ldots i_m}$ is infinite or finite becomes necessary, but not sufficient, for characterizing recurrence or full transience; see the next example.}
\end{itemize}

%In view of the classification results for the first order case recalled at the beginning of this section, it is tempting to raise the question of whether for a higher order chain, we can classify any state $i$ alternatively as recurrent or transient based on whether $\displaystyle \sum_{k=1}^\infty p^{(k)}_{iii_3\ldots i_m}$ is infinite or finite for all $i_3, \ldots, i_m \in S$, in the same manner as the first order case. The next example shows that the answer is negative:

%Denote by $I_3$ the $3$ by $3$ identity matrix. Take a second order chain with three states and transition tensor 
%$${\cal P}(:,:,1)={\cal P}(:,:,3)=I_3, {\cal P}(:,:,2)=\left[\begin{array}{ccc}
%0 & 0 & 1\\
%0 & 1 & 0\\
%1 & 0 & 0
%\end{array}\right].$$
%It is straightforward to check that ${\cal F}=\cal P$. Thus, following the classification scheme in \cite{HX24b}, states $1$ and $3$ are both transient but not fully transient, while state $2$ is recurrent. On the other hand, since ${\cal P}^k={\cal P}$ for all $k \ge 1$, we see  
%$$\sum_{k=1}^\infty {\cal P}^k(:,:,1)=\sum_{k=1}^\infty {\cal P}^k(:,:,3)=\left[\begin{array}{ccc}
%\infty & 0 & 0\\
%0 & \infty & 0\\
%0 & 0 & \infty
%\end{array}\right], \sum_{k=1}^\infty {\cal P}^k(:,:,2)=\left[\begin{array}{ccc}
%0 & 0 & \infty\\
%0 & \infty & 0\\
%\infty & 0 & 0
%\end{array}\right].$$
%Clearly, states 1 and 3 do not satisfy $\displaystyle \sum_{k=1}^\infty p^{(k)}_{iii_3}=\infty$ for all $i_3$. The alternative classification scheme as speculated above, therefore, does not work.

Before proceeding, we need to introduce a few terms. A state $j$ is said to be reachable from state $i$, denoted as $i \rightarrow j$, if there exists $k \ge 0$, which may depend on $j, i, i_3, \ldots, i_m$, such that $p^{(k)}_{jii_3\ldots i_m}>0$ or, equivalently, $f_{jii_3\ldots i_m}>0$, for all $i_3, \ldots, i_m \in S$. In addition, states $i$ and $j$ are said to communicate, denoted by $i \leftrightarrow j$, if $i \rightarrow j$ and $j \rightarrow i$. The $\leftrightarrow$ relationship is an equivalence relationship and thus it partitions the state space $S$ into equivalence classes \cite{HX24b}.

As illustrated in \cite{HX24b}, there are quite a few other differences when we go from first order to higher order chains regarding classification of states. For example, 
\begin{itemize}
\item{A first order chain must have a recurrent state, but not so for a higher order chain. Take a second order chain on three states with transition tensor  
$${\cal P}(:,:,1)=\left[\begin{array}{ccc}
1 & 0 & 0\\
0 & 1/2 & 1/2\\
0 & 1/2 & 1/2
\end{array}\right], ~{\cal P}(:,:,2)=\left[\begin{array}{ccc}
1 & 0 & 1/2\\
0 & 1 & 0\\
0 & 0 & 1/2
\end{array}\right], ~{\cal P}(:,:,3)=\left[\begin{array}{ccc}
0 & 1/2 & 1\\
0 & 1/2 & 0\\
1 & 0 & 0
\end{array}\right].$$
Then,
$${\cal F}(:,:,1)=\left[\begin{array}{ccc}
1 & 1/2 & 3/4\\
0 & 1 & 1\\
0 & 1/2 & 1/2
\end{array}\right], ~{\cal F}(:,:,2)=\left[\begin{array}{ccc}
1 & 0 & 1\\
0 & 1 & 1\\
0 & 0 & 1
\end{array}\right], ~{\cal F}(:,:,3)=\left[\begin{array}{ccc}
3/4 & 1/2 & 1\\
1 & 1/2 & 1\\
1 & 0 & 1
\end{array}\right].$$
For each $1 \le i \le 3$, its ever-reaching probability $f_{iii_3}$ equals $1$ for some $i_3$ and is less than $1$ for some other $i_3$. Consequently, none of the states in this example are recurrent. In addition, we see 
$$\sum_{k=1}^\infty p^{(k)}_{11i_3}=\infty \ \ {\rm and} \ \ \sum_{k=1}^\infty p^{(k)}_{33i_3} < \infty$$
for all $1 \le i_3 \le 3$, implying that whether $\ds \sum_{k=1}^\infty p^{(k)}_{iii_3\ldots i_m}$ is infinite or not for all $i_3, \ldots, i_m \in S$ is only a necessary condition for recurrence or full transience, respectively.
}
\item{In a first order chain, a state $i$ is recurrent provided $j \rightarrow i$ whenever $i \rightarrow j$. Such a characterization does not apply to a higher order chain. Consider a second order chain with two states and transition tensor 
$${\cal P}(:,:,1)=\left[\begin{array}{cc}
1 & 1/2\\
0 & 1/2
\end{array}\right], ~{\cal P}(:,:,2)=\left[\begin{array}{cc}
0 & 1/2\\
1 & 1/2
\end{array}\right].$$
It follows  
$${\cal F}(:,:,1)=\left[\begin{array}{cc}
1 & 1\\
0 & 1
\end{array}\right], ~{\cal F}(:,:,2)=\left[\begin{array}{cc}
1 & 1\\
1 & 1
\end{array}\right].$$
Both states are recurrent since $f_{iii_3}=1$ for all $i, i_3=1, 2$. Meanwhile, $2 \rightarrow 1$ yet $1 \not\rightarrow 2$ since $p^{(k)}_{i_1i_2i_3}>0$ for all $1 \le i_1, i_2 \le 2$ and $k \ge 2$ but $p^{(k)}_{211}=0$ for all $k \ge 1$.}
\item{In a first order chain, both recurrence and transience are class properties, i.e., the states in an equivalence class must be either all recurrent or all transient. When it comes to a higher order chain,  recurrence and transience are not class properties anymore. To see this, we examine a second order chain with three states and transition tensor 
$${\cal P}(:,:,1)=\left[\begin{array}{ccc}
1/2 & 1/3 & 1/2\\
1/2 & 1/2 & 0\\
0 & 1/3 & 1/2
\end{array}\right], ~{\cal P}(:,:,2)=\left[\begin{array}{ccc}
1 & 0 & 1/2\\
0 & 1 & 1/2\\
0 & 0 & 0
\end{array}\right], ~{\cal P}(:,:,3)=\left[\begin{array}{ccc}
0 & 0 & 1/2\\
0 & 0 & 1/2\\
1 & 1 & 0
\end{array}\right].$$
For this chain, we have 
$${\cal F}(:,:,1)=\left[\begin{array}{ccc}
5/6 & 2/3 & 1\\
1 & 1 & 1\\
1/2 & 1/2 & 1
\end{array}\right], ~{\cal F}(:,:,2)=\left[\begin{array}{ccc}
1 & 0 & 1\\
1 & 1 & 1\\
1/2 & 0 & 1
\end{array}\right], ~{\cal F}(:,:,3)=\left[\begin{array}{ccc}
1 & 1 & 1\\
1 & 1 & 1\\
1 & 1 & 1
\end{array}\right].$$
By checking its ever-reaching probabilities $f_{iii_3}$, state $1$ is transient but not fully transient, and state $3$ is recurrent. These states, however, are in the same equivalence class, i.e., $1 \leftrightarrow 3$.}
\end{itemize}
Such differences constitute yet another roadblock to studying a higher order chain via its reduced first order chain. Thus, as far as classification of states is concerned, the answer to Question \ref{ques} is again negative.

The following was conjectured in \cite{HX24b}. To end this section, we shall prove it here.
\begin{2}
\label{conj2}
Let $K \subset S$ be an equivalence class. Then, it is impossible to have $i, j \in K$ such that $i$ is recurrent while $j$ is fully transient.
\end{2}
\bp
Let $i, j \in K$, $i \ne j$, with $j$ being recurrent. Then, it suffices to show that $i$ cannot be fully transient.

Let us fix $i_3, \ldots, i_m$. Since $i \rightarrow j$ and $i \ne j$, there exists $\alpha \ge 1$ such that $p^{(\alpha)}_{jii_3\ldots i_m}>0$. Denote the event ``the chain takes $\alpha$ steps to go from $ii_3\ldots i_m$ to $j$'' by $A$.

Next, denote the event ``starting from $jj_3\ldots j_m$, the chain takes $\beta$ steps to return to $j$'' by $B$, where $j_3, \ldots, j_m$ are given by:

\begin{itemize}
\item{when $\alpha=1$: $j_3=i, j_4=i_3, \ldots, j_m=i_{m-1}$}
\item{when $\alpha=2$: $j_4=i, j_5=i_3, \ldots, j_m=i_{m-2}$, any $j_3$}
\item{$\cdots$}
\item{when $\alpha=m-2$: $j_m=i$, any $j_3,\ldots,j_{m-1}$}
\item{$\cdots$}
\end{itemize}

Fix $j_3, \ldots, j_m$ from now on.

Next, denote the event ``the chain takes $\gamma$ steps to go from $jk_3\ldots k_m$ to $i$'' by $C$, where $k_3, \ldots, k_m$ are determined in the same way as $j_3, \ldots, j_m$ above.

From this point on, let us assume $\beta \ge m-1$. Thus, we can fix $k_3, \ldots, k_m$. Accordingly, by $j \rightarrow i$ and $i \ne j$, there exists $\gamma \ge 1$ such that $p^{(\gamma)}_{ijk_3\ldots k_m}>0$.

Finally, we observe $$\begin{array}{rcl}
p^{(\gamma+\beta+\alpha)}_{iii_3\ldots i_m} & \ge & \Pr(ABC)\\
 & = & \Pr(A)\Pr(B|A)\Pr(C|AB)\\
 & = & p^{(\alpha)}_{jii_3\ldots i_m}p^{(\beta)}_{jjj_3\ldots j_m}p^{(\gamma)}_{ijk_3\ldots k_m}.
\end{array}$$
It follows that 
$$\sum_{\beta=m-1}^\infty p^{(\gamma+\beta+\alpha)}_{iii_3\ldots i_m} \ge p^{(\alpha)}_{jii_3\ldots i_m}p^{(\gamma)}_{ijk_3\ldots k_m}\sum_{\beta=m-1}^\infty p^{(\beta)}_{jjj_3\ldots j_m}=\infty$$
since $j$ is recurrent. The above clearly leads to $$\sum_{\beta=1}^\infty p^{(\beta)}_{iii_3\ldots i_m}=\infty.$$
The proof is now complete.
\ep

\section{Mean First Passage Times}
\label{mfpt}
\setcounter{equation}{0}

For convenience, let us use a framework that includes both higher order and first order chains.

Recall $\eta_{i_1i_2\ldots i_m}$ in (\ref{eta}). For an $(m-1)$th order, $m \ge 2$, chain with state space $S=\{1, 2, \ldots, n\}$, the mean first passage time from states $(i_2, \ldots, i_m)$ to $i_1$ is defined by  
$$\mu_{i_1i_2\ldots i_m}={\rm E}(\eta_{i_1i_2\ldots i_m})=\sum_{k=1}^\infty k\Pr(\eta_{i_1i_2\ldots i_m}=k).$$
Accordingly, the mean first passage time tensor is the $m$th order, $n$ dimensional tensor $\mu=[\mu_{i_1i_2\ldots i_m}]$.

As shown in \cite{HX24a}, for an ergodic higher order chain of order $m-1$ whose $m$th order, $n$ dimensional transition tensor is $\cal P$, its mean first passage times are finite and are uniquely determined by the following tensor equation:
\be
\label{mfp_eqn}
\mu={\cal E}+(\mu-\mu_d)\boxtimes {\cal P},
\ee
where $\cal E$ is the $m$th order, $n$ dimensional tensor of all ones. Furthermore, as a linear system, (\ref{mfp_eqn}) is nonsingular if and only if the chain is ergodic.

Comparing with an ergodic first order chain, the mean first passage time matrix $M$ satisfies \cite{KS}
\be
\label{mfp_eqn_1}
M=E+(M-M_d)P,
\ee
where $E$ is the matrix of all ones and $M_d=[m^{(d)}_{ij}]$ is such that 
$$m^{(d)}_{ij}=\left\{\begin{array}{cl}
m_{ii}, & i=j;\\
0, & {\rm otherwise}.
\end{array}\right.$$
Clearly, (\ref{mfp_eqn_1}) is a special case of (\ref{mfp_eqn}). Because of the resemblance between (\ref{mfp_eqn}) and (\ref{mfp_eqn_1}), and in the spirit of Question \ref{ques}, it is an intriguing question whether (\ref{mfp_eqn}) can be derived from (\ref{mfp_eqn_1}) using the reduced first order chain. The answer, however, turns out to be negative, as illustrated by the example below.

Take a second order chain with three states and $${\cal P}(:,:,1)={\cal P}(:,:,2)={\cal P}(:,:,3)=\left[\begin{array}{ccc}
1/3 & 1/3 & 1/3\\
1/3 & 1/3 & 1/3\\
1/3 & 1/3 & 1/3
\end{array}\right].$$
By solving (\ref{mfp_eqn}), $$\mu(:,:,1)=\mu(:,:,2)=\mu(:,:,3)=\left[\begin{array}{ccc}
3 & 3 & 3\\
3 & 3 & 3\\
3 & 3 & 3
\end{array}\right].$$
On the other hand, observing that the regularity of this second order chain carries over to its reduced first order chain since $Q^2=(1/9)E>0$, where $E$ is the matrix of all ones, and using (\ref{mfp_eqn_1}) for the reduced first order chain, we obtain  
$$M=\left[\begin{array}{ccccccccc}
    9 &  12 &  12 &   9 &  12 &  12 &   9 &  12 &  12\\
    6 &   9 &   9 &   6 &   9 &   9 &   6 &   9 &   9\\
    6 &   9 &   9 &   6 &   9 &   9 &   6 &   9 &   9\\
    9 &   6 &   9 &   9 &   6 &   9 &   9 &   6 &   9\\
   12 &   9 &  12 &  12 &   9 &  12 &  12 &   9 &  12\\
    9 &   6 &   9 &   9 &   6 &   9 &   9 &   6 &   9\\
    9 &   9 &   6 &   9 &   9 &   6 &   9 &   9 &   6\\
    9 &   9 &   6 &   9 &   9 &   6 &   9 &   9 &   6\\
   12 &  12 &   9 &  12 &  12 &   9 &  12 &  12 &   9
\end{array}\right].$$
Obviously, in this example, $\mu$ does not follow directly from $M$. If we look at $m_{11}$, for instance, it represents the mean first passage time from state $11$, in multi-index form, to itself. This quantity, however, is different from $\mu_{111}$. In other words, while addressing Question \ref{ques}, other difficulties may still arise even if ergodicity or regularity does carry forward from a higher order chain to its reduced first order chain.

One way of working around the discrepancy as shown above between the mean first passage times of a higher order chain and those of its reduced first order chain is to resort to certain matricization of the transition tensor \cite{HX24a}. This approach is not as convenient as the direct tensor treatment leading to (\ref{mfp_eqn}). Besides, its result can be easily obtained by matricizing the mean first passage time tensor from (\ref{mfp_eqn}).

\section{Limiting Probability Distributions}
\label{limi}
\setcounter{equation}{0}

Having presented a range of differences between higher order and first order chains in previous sections, we now come to a point where these two appear to converge to a certain level. Again, we use here a framework that includes both higher order and first order chains.

With $m \ge 2$, let $X=\{X_t : t=1, 2, \ldots\}$ be an $(m-1)$th order chain on state space $S=\{1, 2, \ldots, n\}$ whose transition tensor is ${\cal P}=[p_{i_1i_2\ldots i_m}]$. The probability distribution of $X$ at $t$ is given by 
$$x_t=[\Pr(X_t=1) ~\Pr(X_t=2) ~\ldots ~\Pr(X_t=n)]^T.$$
If there exists an $n$-vector $\pi$, which is independent of the initial probability distributions $x_1, x_2, \ldots, x_{m-1}$, such that 
$$\lim_{t \rightarrow \infty}x_t=\pi,$$
then $\pi$ is called the limiting probability distribution of the chain $X$\cite{HX25}. It is clear that $\pi$ satisfies $\pi \ge 0$ and $\ds \sum_{i \in S}\pi_i=1$.

On the other hand, let $Y$ be the reduced first order chain associated with $X$, whose transition matrix is given by $Q$. Then, an $N$-vector $\xi$ is said to be a stationary distribution of both the chains $X$ and $Y$ if $\xi$ satisfies $\xi \ge 0$, $\ds \sum_{i_1i_2\ldots i_{m-1} \in T}\xi_{i_1i_2\ldots i_{m-1}}=1$, and $Q\xi=\xi$. The subscript for the entries of $\xi$ is written here in multi-index form. Consequently, $Q\xi=\xi$ translates to  
$$\sum_{j \in S}p_{i_1i_2\ldots i_{m-1}j}\xi_{i_2\ldots i_{m-1}j}=\xi_{i_1i_2\ldots i_{m-1}}, ~i_1, i_2, \ldots, i_{m-1} \in S,$$
which is precisely how a stationary distribution is equivalently defined in \cite{GLY} when $\xi$ is treated as an $(m-1)$th order, $n$ dimensional tensor. Obviously, as a vector, $\xi$ is just a normalized right eigenvector of $Q$ corresponding to the dominant eigenvalue $\lambda=1$, sometimes also called a Perron vector \cite{HJ}.

It is shown in \cite{HX25} that if the chain $X$ is regular, then its limiting probability distribution $\pi$ is guaranteed to exist. In this case, $\pi>0$ and 
$$\lim_{k \rightarrow \infty}{\cal P}^k=\pi \otimes \underbrace{e \otimes \cdots \otimes e}_{m-1},$$
where, and in what follows, $e$ is the $n$-vector of all ones and $\otimes$ is the outer product. Without the regularity condition on $X$, $\pi$ may not exist and therefore such a condition is tight. 

Regarding a stationary distribution of the chains $X$ and $Y$, its existence is guaranteed following the celebrated Perron-Frobenius theorem, yet such a distribution may not be unique. If, in addition, $Y$ is a regular chain, then its stationary distribution $\xi$ is unique and $\xi > 0$.

The preceding results regarding the limiting probability distribution and stationary distribution continue to hold when the chain $X$ is first order and regular \cite{KS, Hun}. In this scenario, the limiting probability distribution is the same as the stationary distribution. To be more specific, on letting $P$ be the transition matrix, then there exists a unique $n$-vector $\pi > 0$ with $\ds \sum_{i \in S}\pi_i=1$ such that $\ds \lim_{k \rightarrow \infty}P^k=\pi e^T$ and, at the same time, $P\pi=\pi$ \cite{KS}. This appears to be the cause of why the nomenclatures ``limiting probability distribution'' and ``stationary distribution'' are at times used interchangeably in works on higher order chains. Our foregoing description of these concepts clarifies the difference between them in the higher order case.

Despite the differences, there is a principal connection between them too. As shown in \cite{HX25}, when the chain $X$ is higher order and regular, its limiting probability distribution $\pi$ can be determined via any stationary distribution $\xi$ through 
\be
\label{pixi}
\pi=P^{(0)}\xi,
\ee
where $P^{(0)}$ stands for the mode-$1$ matricization of the $m$th order, $n$ dimensional identity tensor. This result does not require the chain $Y$ to be regular.

To illustrate the usage of (\ref{pixi}), let us take the same regular second order chain as in Section \ref{irre} with four states and transition tensor  
$${\cal P}(:,:,1)=\left[\ba{cccc}
1/2 & 0 & 0 & 0\\
1/2 & 0 & 1 & 0\\
0 & 1 & 0 & 1\\
0 & 0 & 0 & 0
\ea\right], ~{\cal P}(:,:,2)=\left[\ba{cccc}
0 & 0 & 1/2 & 1\\
0 & 1/2 & 0 & 0\\
1/2 & 1/2 & 0 & 0\\
1/2 & 0 & 1/2 & 0
\ea\right],$$
$${\cal P}(:,:,3)=\left[\ba{cccc}
0 & 1 & 0 & 1\\
1 & 0 & 1/2 & 0\\
0 & 0 & 1/2 & 0\\
0 & 0 & 0 & 0
\ea\right], ~{\cal P}(:,:,4)=\left[\ba{cccc}
0 & 0 & 0 & 0\\
1 & 1 & 1 & 0\\
0 & 0 & 0 & 1/2\\
0 & 0 & 0 & 1/2
\ea\right].$$
The associated reduced first order chain is reducible, which can also be seen from the fact that the dominant eigenvalue $\lambda=1$ of its transition matrix $Q$ has multiplicity $2$. This, of course, implies that the reduced first order chain is not regular. It is easy to verify that $$\displaystyle z=\frac{1}{7}[0 ~~0 ~~1 ~~1 ~~2 ~~0 ~~0 ~~0 ~~0 ~~2 ~~0 ~~0 ~~0 ~~0 ~~1 ~~0]^T$$ is a normalized right eigenvector corresponding to that dominant eigenvalue above. Hence, the limiting probability distribution of the second order chain is
$$\pi=P^{(0)}z=\frac{1}{7}[2 ~~2 ~~2 ~~1]^T.$$

Regarding Question \ref{ques}, when it comes to the limiting probability distribution problem of a higher order chain, (\ref{pixi}) represents its close tie with a stationary distribution of the corresponding reduced first order chain. This, along with the $k$-step transition probabilities in Section \ref{prob}, is another case in which Question \ref{ques} yields a positive answer.

\section{Conclusions}
\label{concl}
\setcounter{equation}{0}

In this article, we provide a relatively comprehensive answer to Question \ref{ques} in order to address a still prevailing perception, or ``rooted prejudice'' in the words of Marius Iosifescu, concerning higher order chains. Our results show that there can be significant differences between a higher order chain and its reduced first order chain in most aspects such as ergodicity, regularity, ever-reaching probabilities, mean first passage times, and classification of states. Owing to such differences, the answer to Question \ref{ques} is negative in general. In the meantime, the closest, and also probably the most useful, tie between a higher order chain and its reduced first order chain is seen in the limiting probability distribution problem. Overall, it remains a very interesting topic to further explore whether a reduced first order chain may, or may not, play a useful role in studying other key quantities for its originating higher order chain.


\begin{thebibliography}{}

\bibitem{Bae}
S. Baena-Mirabete, P. Puig, Parsimonious higher order Markov models for rating transitions, {\it Journal of the Royal Statistical Society A: Statistics in Society} 181: 107--131, 2018.

\bibitem{Bur}
D. Burks, R. Azad, Higher-order Markov models for metagenomic sequence classification, {\it Bioinformatics} 36: 4130--4136, 2020.

\bibitem{CZ}
K. Chang, T. Zhang, On the uniqueness and non-uniqueness of the positive $Z$-eigenvector for transition probability tensors, {\it Journal of Mathematical Analysis \& Applications} 408: 525--540, 2013.

\bibitem{CPZ}
J. Culp, K. Pearson, T. Zhang, On the uniqueness of the $Z_1$-eigenvector of transition probability tensors, {\it Linear \& Multilinear Algebra} 65: 891--896, 2017.

\bibitem{Doo}
J. Doob, {\it Stochastic Processes}, Wiley Classic Library Edition, Wiley, 1990.

\bibitem{Fle}
G. Flett, N. Kelly, An occupant-differentiated, higher-order Markov chain method for prediction of domestic occupancy, {\it Energy \& Buildings} 125: 219--230, 2016.

\bibitem{G}
B. Geiger, A sufficient condition for a unique invariant distribution of a higher-order Markov chain, {\it Statistics \& Probability Letters} 130: 49--56, 2017.

\bibitem{GLY}
D. Gleich, L. Lim, Y. Yu, Multilinear pagerank, {\it SIAM Journal on Matrix Analysis \& Applications} 36: 1507--1541, 2015.

\bibitem{Ho}
L. Ho, J. Rajapakse, Splice site detection with a higher-order Markov model implemented on a neural network, {\it Genome Informatics} 14: 64--72, 2003.

\bibitem{HQ}
S. Hu, L. Qi, Convergence of a second order Markov chain, {\it Applied Mathematics \& Computation} 241: 183--192, 2014.

\bibitem{Hun}
J. Hunter, {\it Mathematical Techniques of Applied Probability, Volume I: Discrete Time Models: Basic Theory}, Academic Press, 1983.

\bibitem{HWX}
L. Han, K. Wang, J. Xu, Higher order ergodic Markov chains and first passage times, {\it Linear \& Multilinear Algebra} 70: 6772--6779, 2022.

\bibitem{HX24a}
L. Han, J. Xu, Ever-reaching probabilities and mean first passage times of higher order ergodic Markov chains, {\it Linear \& Multilinear Algebra} 72: 59--75, 2024.

\bibitem{HX24b}
L. Han, J. Xu, On classification of states in higher order Markov chains, {\it Linear Algebra \& Its Applications} 685: 24--45, 2024.

\bibitem{HX25}
L. Han, J. Xu, On limiting probability distributions of higher order Markov chains, under review, \url{https://doi.org/10.48550/arXiv.2506.08874}.

\bibitem{HJ}
R. Horn, C. Johnson, {\it Matrix Analysis}, Cambridge University Press, 1990.

\bibitem{Ios}
M. Iosifescu, {\it Finite Markov Processes \& Their Applications}, Dover Publications, 2007.

\bibitem{Isl}
M. Islam, R. Chowdhury, A higher order Markov model for analyzing covariate dependence, {\it Applied Mathematical Modelling} 30: 477--488, 2006.

\bibitem{KS}
J. Kemeny, J. Snell, {\it Finite Markov Chains}, Springer-Verlag, 1960.

\bibitem{KB}
T. Kolda, B. Bader, Tensor decompositions and applications, {\it SIAM Review} 51: 455--500, 2009.

\bibitem{Kwa}
J. Kwak, C. Lee, D. Eun, A high-order Markov-chain-based scheduling algorithm for low delay in CSMA networks, {\it IEEE/ACM Transactions on Networking} 24: 2278--2290, 2015.

\bibitem{Lan}
J. Lan, X. Li, V. Jilkov, C. Mu, Second-order Markov chain based multiple-model algorithm for maneuvering target tracking, {\it IEEE Transactions on Aerospace \& Electronic Systems} 49: 3--19, 2013.

\bibitem{LZ}
C. Li, S. Zhang, Stationary probability vectors of higher-order Markov chains, {\it Linear Algebra \& Its Applications} 473: 114--125, 2016.

\bibitem{LN}
W. Li, M. Ng, On the limiting probability distribution of a transition probability tensor, {\it Linear Algebra \& Its Applications} 62: 362--385, 2014.

\bibitem{LLZ}
Z. Liu, Y. Luo, Y. Zhu, State estimation for linear dynamic system with multiple-step random delays using higher-order Markov chain, {\it IEEE Access} 8: 76218--76227, 2020.

\bibitem{MSL}
C. Martin, R. Shafer, B. Larue, An order-$p$ tensor factorization with applications in imaging, {\it SIAM Journal on Scientific Computing} 35: 474--490, 2013.

\bibitem{M}
N. Masseran, Markov chain model for the stochastic behaviors of wind-direction data, {\it Energy Conversion \& Management} 92: 266--274, 2015.

\bibitem{San}
M. Sanjari, H. Gooi, Probabilistic forecast of PV power generation based on higher order Markov chain, {\it IEEE Transactions on Power Systems} 32: 2942--2952, 2016.

\bibitem{WC}
S. Wu, M. Chu, Markov chains with memory, tensor formulation, and the dynamics of power iteration, {\it Applied Mathematics \& Computation} 303: 226--239, 2017.

\bibitem{Xio}
H. Xiong, R. Mamon, A higher-order Markov chain-modulated model for electricity spot-price dynamics, {\it Applied Energy} 233-234: 495--515, 2019.

\bibitem{X25}
J. Xu, HOMC: a MATLAB package for higher order Markov chains, under review, \url{https://doi.org/10.48550/arXiv.2510.02664}.

\bibitem{YJK}
Y. Yang, H. Jang, B. Kim, A hybrid recommender system for sequential recommendation: combining similarity models with Markov chains, {\it IEEE Access} 8: 190136--190146, 2020.

\end{thebibliography}
\end{document}